\documentclass{amsart} 
\usepackage{amscd,amssymb,latexsym,subfigure,verbatim,hyperref,epsfig}
\usepackage[graph,frame,poly,arc]{xy}
\usepackage{calc}
\usepackage{amsmath}
\usepackage{graphicx}

\begin{document}

\newcommand{\mmbox}[1]{\mbox{${#1}$}}
\newcommand{\proj}[1]{\mmbox{{\mathbb P}^{#1}}}
\newcommand{\affine}[1]{\mmbox{{\mathbb A}^{#1}}}
\newcommand{\Ann}[1]{\mmbox{{\rm Ann}({#1})}}
\newcommand{\caps}[3]{\mmbox{{#1}_{#2} \cap \ldots \cap {#1}_{#3}}}
\newcommand{\N}{{\mathbb N}}
\newcommand{\Z}{{\mathbb Z}}
\newcommand{\R}{{\mathbb R}}
\newcommand{\K}{{\mathbb K}}
\newcommand{\p}{{\mathbb P}}
\newcommand{\A}{{\mathcal A}}
\newcommand{\CC}{{\mathcal C}}
\newcommand{\C}{{\mathbb C}}

\newcommand{\Tor}{\mathop{\rm Tor}\nolimits}
\newcommand{\Ext}{\mathop{\rm Ext}\nolimits}
\newcommand{\Hom}{\mathop{\rm Hom}\nolimits}
\newcommand{\im}{\mathop{\rm Im}\nolimits}
\newcommand{\reg}{\mathop{\rm reg}\nolimits}
\newcommand{\syz}{\mathop{\rm syz}\nolimits}
\newcommand{\pdim}{\mathop{\rm pdim}\nolimits}
\newcommand{\Sing}{\mathop{\rm Sing}\nolimits}
\newcommand{\rank}{\mathop{\rm rank}\nolimits}
\newcommand{\supp}{\mathop{\rm supp}\nolimits}
\newcommand{\arrow}[1]{\stackrel{#1}{\longrightarrow}}
\newcommand{\CB}{Cayley-Bacharach}
\newcommand{\coker}{\mathop{\rm coker}\nolimits}
\sloppy
\newtheorem{defn0}{Definition}[section]
\newtheorem{prop0}[defn0]{Proposition}
\newtheorem{conj0}[defn0]{Conjecture}
\newtheorem{thm0}[defn0]{Theorem}
\newtheorem{lem0}[defn0]{Lemma}
\newtheorem{corollary0}[defn0]{Corollary}
\newtheorem{example0}[defn0]{Example}

\newenvironment{defn}{\begin{defn0}}{\end{defn0}}
\newenvironment{prop}{\begin{prop0}}{\end{prop0}}
\newenvironment{conj}{\begin{conj0}}{\end{conj0}}
\newenvironment{thm}{\begin{thm0}}{\end{thm0}}
\newenvironment{lem}{\begin{lem0}}{\end{lem0}}
\newenvironment{cor}{\begin{corollary0}}{\end{corollary0}}
\newenvironment{exm}{\begin{example0}\rm}{\end{example0}}

\newcommand{\defref}[1]{Definition~\ref{#1}}
\newcommand{\propref}[1]{Proposition~\ref{#1}}
\newcommand{\thmref}[1]{Theorem~\ref{#1}}
\newcommand{\lemref}[1]{Lemma~\ref{#1}}
\newcommand{\corref}[1]{Corollary~\ref{#1}}
\newcommand{\exref}[1]{Example~\ref{#1}}
\newcommand{\secref}[1]{Section~\ref{#1}}

\newcommand{\poina}{\pi({\mathcal A}, t)}
\newcommand{\poinc}{\pi({\mathcal C}, t)}

\newcommand{\std}{Gr\"{o}bner}
\newcommand{\jq}{J_{Q}}



\title {Freeness of Conic-Line arrangements in $\mathbb{P}^2$}

\author{Hal Schenck}
\thanks{Schenck supported by NSF 03--11142, 07--07667, NSA 904-03-1-0006}\address{Schenck: Mathematics Department \\ Texas A\&M University \\
  College Station \\ TX 77843-3368\\ USA, and Mathematics Department \\ University of Illinois \\
   Urbana \\ IL 61801\\  }
\email{schenck@math.uiuc.edu}

\author{\c Stefan O. Toh\v aneanu}
\address{Tohaneanu: Department of Mathematical Sciences \\ University of Cincinnati \\ Cincinnati \\ OH 45221-0025\\ }
\email{stefan.tohaneanu@uc.edu}

\subjclass[2000]{Primary 52B30; Secondary 14J60} \keywords{line
  arrangement, curve arrangement, module of derivations.}

\begin{abstract}
\noindent Let ${\mathcal C}= \bigcup_{i=1}^n C_i \subseteq
\mathbb{P}^2$
be a collection of smooth rational plane curves. We prove that
the addition-deletion operation used in the study of hyperplane
arrangements has an extension which works for a large class of
arrangements of smooth rational curves, giving an inductive tool
for understanding the freeness of the module
$\Omega^1({\mathcal C})$ of logarithmic differential forms
with pole along ${\mathcal C}$. We also show that the
analog of Terao's conjecture (freeness of  $\Omega^1({\mathcal C})$
is combinatorially determined if ${\mathcal C}$ is a union of lines)
is false in this setting.
\end{abstract}
\maketitle


\section{Introduction}\label{sec:intro}
One of the fundamental objects associated to a hyperplane arrangement
${\mathcal A} \subseteq \p_{\K}(V)$ is the module $\Omega^1({\mathcal
  A})$
of logarithmic one-forms with pole along the arrangement or (dually) the module  $D({\mathcal A})$ of
derivations tangent to the arrangement. Both are graded $S = Sym(V^*)$ modules; $D(\A) \subseteq Der_{\K}(S)$ is
defined via:
\begin{defn}
 $D({\mathcal A}) = \{ \theta | \theta(l_i) \in \langle l_i \rangle \}$
  for all $l_i$ such that $V(l_i) \in \A$.
\end{defn}
Over a field of characteristic zero,
$D({\mathcal A})\simeq E \oplus D_0({\mathcal A})$, where
$E$ is the Euler derivation and $D_0({\mathcal A})$
corresponds to the module of syzygies on the Jacobian ideal
of the defining polynomial of ${\mathcal A}$. When $\K = \C$ or $\R$,
an elegant theorem of Terao relates the freeness of the
module $D({\mathcal A})$ to the Poincar\'e polynomial of
$V \setminus \A$. In this note, we restrict to $\mathbb{P}^2$, but
broaden the class of curves which make up the arrangement.
In particular, suppose
\[
\mathcal{C} = \bigcup_{i=1}^nC_i,
\]
where each $C_i$ is a smooth rational plane curve; call such a collection
a conic-line (CL) arrangement.

\begin{exm}\label{firstex}
For the CL arrangement below, $D({\mathcal C}) \simeq
S(-1) \oplus S(-2)\oplus S(-5)$.
\begin{figure}[h]
\begin{center}
\epsfig{file=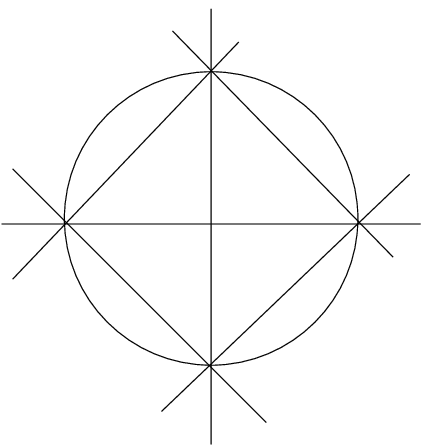,height=1.0in,width=1.0in}
\end{center}
\end{figure}
\end{exm}
\subsection{Line arrangements}
We begin with some facts about hyperplane arrangements; for more information  see Orlik and Terao \cite{ot}.
\noindent A hyperplane arrangement ${\mathcal A}$ is a finite collection of codimension one linear subspaces of
a fixed vector space V. ${\mathcal A}$ is {\it central} if each hyperplane contains the origin {\bf 0} of V. The
intersection lattice $L_{\mathcal A}$ of ${\mathcal A}$ consists of the intersections of the elements of
${\mathcal A}$; the rank of $x \in L_{\mathcal A}$ is simply the codimension of $x$. V is the lattice element
$\hat{0}$; the rank one elements are the hyperplanes themselves. ${\mathcal A}$ is called {\it essential} if
rank $L_{\mathcal A} =$ dim $V$.
\begin{defn}
The M\"{o}bius function $\mu$ : $L_{\mathcal A} \longrightarrow \Z$ is defined by $$\begin{array}{*{3}c}
\mu(\hat{0}) & = & 1\\
\mu(t) & = & -\sum\limits_{s < t}\mu(s) \mbox{, if } \hat{0}< t
\end{array}$$
\end{defn}
We now restrict to the case that $V$ is complex. A foundational result
is that the Poincar\'e polynomial of $X = V \setminus {\mathcal A}$
is purely combinatorial; in particular
\[
P(X,t) = \sum\limits_{x \in L_{\mathcal A}}\mu(x) \cdot (-t)^{\text{rank}(x)}.
\]
An arrangement $\mathcal{A}$ is {\em free} if $D({\mathcal A}) \simeq
\oplus S(-a_i)$; the $a_i$ are called the {\em exponents} of ${\mathcal A}$.
Terao's famous theorem \cite{t} is that if $D({\mathcal A}) \simeq
\oplus S(-a_i)$, then $P(X,t)= \prod(1+a_it)$.
If $\A \subseteq \C^3$ is central, then $\A$ also defines a set of
lines in $\p^2$, and obviously $X = \C^3
\setminus \A \simeq \C^{*} \times \widetilde{X}$, where
$\widetilde{X}$ is the complement of the corresponding
arrangement of lines in $\mathbb{P}^2$. Hence
\[
P(\widetilde{X},t)= 1+(n-1)t+(\!\!\!\!\!\sum\limits_{\stackrel{x \in
    L_{\mathcal A}}{\text{rank}(x)=2}}\!\!\!\! \mu(x) -n+1) t^2
\]
It follows from Terao's theorem that if
$D_0({\mathcal A}) \simeq S(-a)\oplus S(-b)$, then
$P(\widetilde{X},t) = (1 + at)(1+bt)$. This can be generalized
\cite{sch} to line arrangements which are not free, using the
Chern polynomial. The motivating question of this paper is:
{\em what happens if the arrangement of lines is replaced with a CL arrangement?}

\subsection{Rational curve arrangements}
In \cite{cg}, Cogolludo-Agust\'{\i}n studies the complement of an arrangement
of rational curves in $\p^2$, where the individual curves can have
singularities, and can meet non-transversally. The main result is
that the cohomology ring of the complement to a rational curve
arrangement is generated by logarithmic $1$ and $2$-forms and its
structure depends on a finite number of invariants of the curve.
One fact is that if $\widetilde{X}$ is the complement of an arrangement
of $n$ irreducible curves in $\p^2$, then
\[
h^1(\widetilde{X},\C) = n-1
\]
\[
h^2(\widetilde{X},\C) = 1 + \sum\limits_{p \in \Sing(\CC)}(r_p -1) - \sum\limits_{1}^n (\chi(\hat{C_i})-1),
\]
where $r_p$ is the number of branches passing thru $p$, and
$\hat{C_i}$ is the normalization of $C_i$. Since we are assuming
that all the $C_i$ are smooth and rational, we have that
\[
h^2(\widetilde{X},\C) = \sum\limits_{p \in L_2(\CC)}(r_p-1) - |\CC| +1,
\]
where the intersection poset $L(\CC)$ is defined precisely as for a
linear arrangement (typically, $L(\CC)$ is only a poset, not a lattice).
\subsection{Milnor and Tjurina numbers}
A crucial distinction between line and curve arrangements, even in our
simple setting, is the difference between the Milnor and Tjurina
numbers at a singularity. Let $C=V(f)$ be a reduced
(but not necessarily irreducible) curve in $\mathbb{C}^2$, let
$(0,0) \in C$, and let $\mathbb{C}\{x,y\}$ denote the ring of
convergent power series.
\begin{defn}
The Milnor number of $C$ at $(0,0)$ is
\[
\mu_{(0,0)}(C) = \dim_\mathbb{C}\mathbb{C}\{x,y\}/\langle\frac{\partial f}{\partial
x}\mbox{, }\frac{\partial f}{\partial y}\rangle.
\]
\end{defn}
\noindent To define $\mu_p$ for an arbitrary point $p$, we translate so
that $p$ is the origin.
\begin{defn}
The Tjurina number of $C$ at $(0,0)$ is
\[
\tau_{(0,0)}(C) = \dim_{\C} \C \{x,y\}/\langle\frac{\partial f}{\partial
x}\mbox{, }\frac{\partial f}{\partial y}\mbox{, }f\rangle.
\]
\end{defn}
\begin{defn}
A singularity is quasihomogeneous iff there exists a holomorphic
change of variables so the defining equation becomes weighted
homogeneous; $f(x,y) = \sum c_{ij}x^{i} y^{j}$ is
weighted homogeneous if there exist rational numbers $\alpha, \beta$
such that $\sum c_{ij}x^{i \cdot \alpha} y^{j \cdot \beta}$ is homogeneous.
\end{defn}
In \cite{r}, Reiffen proved that if $f(x,y)$ is a convergent power series
with isolated singularity at the origin, then $f(x,y)$ is in the ideal
generated by the partial derivatives if and only if $f$ is
quasihomogeneous (see \cite{s} for a generalization).

As noted earlier, for a line arrangement with defining
polynomial $F$, $D_0(\mathcal{A})$ consists of the
syzygies on the Jacobian ideal $J_F$ of $F$. If $V(F) \subseteq
\mathbb{P}^2$ is a reduced curve, then after a change of coordinates,
we may assume that $V(F)$ has no singularities on the
line $z=0$. Dehomogenizing so that $f(x,y)= F(X,Y,1)$ yields:
\[
\text{deg}(J_F) = \text{dim}_{\mathbb{C}}\mathbb{C}[x,y]/\langle
\frac{\partial f}{\partial x},\frac{\partial f}{\partial y}, f \rangle
= \!\!\!\!\sum\limits_{p \in \text{Sing}(C)}\!\!\!\!\tau_p(f).
\]
It follows that if all the singular points are quasihomogeneous then
\[
\text{deg}(J_F) = \!\!\!\!\sum\limits_{p \in \text{Sing}(V(f))}\!\!\!\!\mu_p(f).
\]
For a line arrangement, the singularities are always quasihomogeneous,
but this is not the case for CL arrangements:

\begin{exm}
Let $\mathcal C = V(xy(x-y)(x-2y)(x^2-xz+y^2-yz))$ be as below:
\begin{figure}[h]
\begin{center}
\includegraphics{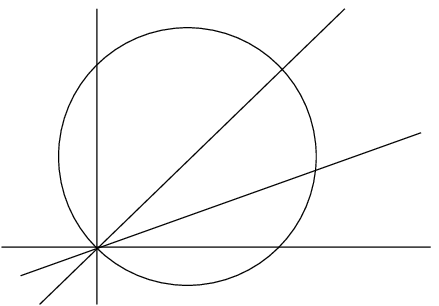}
\end{center}

\end{figure}

\noindent $\mathcal C$ has five singular points, all ordinary.
When $p$ is an ordinary singularity and $C$ has $n$ distinct
branches at $p$, then $\mu_p(C)=(n-1)^2$, so the sum of the
Milnor numbers is 20. However, $\deg(J)=19$;
at $(0:0:1)$ we have $\mu=16$ but $\tau=15$.
\end{exm}

\subsection{Criteria for freeness}\label{sec:two}
The first criterion for the freeness of $D(A)$ is
\begin{prop}[Saito, \cite{slog}]
$\mathcal{A}$ is free exactly when there exist $n+1$ elements
\[
\theta_i = \sum\limits_{j=0}^n f_{ij}\frac{\partial}{\partial x_j} \in D(\mathcal{A})
\]
such that the determinant of the matrix $[f_{ij}]$ is a nonzero
constant multiple of the defining polynomial of $\mathcal{A}$.
\end{prop}
Saito's criterion holds for an arrangement of reduced hypersurfaces
 $\mathcal{C} \subseteq \mathbb{P}^n$;
let $\mathcal{C} = V(F)$ where $F=f_1\cdots f_d$, and $gcd(f_i,f_j)=1$.
By induction,
\[
\theta(f_1\cdots f_d)=f_1\theta(f_2\cdots f_d)+f_2\cdots
f_d\theta(f_1) \in \langle f_1\cdots f_d\rangle,
\]
so we have
\[
D(\mathcal C)= \{ \theta \in Der_{\mathbb K}(S) \mid \theta(f_i) \in
\langle f_i \rangle, i=1,\ldots,d\} = \{\theta\in Der_{\mathbb K}(S)
\mid \theta(F) \in \langle
F\rangle \}.
\]
Any arrangement of (reduced) hypersurfaces will have a singular locus of
codimension two. As for a linear arrangement, $D(\mathcal{C}) \simeq E
\oplus D_0(\mathcal{C})$, with $D_0(\mathcal{C}) = \syz(J_F)$, so
freeness is equivalent to
$\text{pdim}(S/J_F)=2$ (so also equivalent to $J_F$ Cohen-Macaulay).
By the Hilbert-Burch theorem (\cite{e}),
any codimension two Cohen-Macaulay ideal $I$ with $m+1$ generators
is generated by the maximal minors of an $m \times m+1$ matrix $M$,
whose columns generate the module of first syzygies on $I$.
So when $I = J_F$, appending a column vector
$[x_0,\ldots,x_n]$ to $M$ and taking the determinant yields a
multiple of $F$, by Euler's formula. Saito's criterion is
most useful when an explicit set of candidates for the generating set of
$\syz(J_F)$ is known.

There are two other fundamental tools that can be used to prove that
a line arrangement is free. The first method is based on an inductive
operation known as deletion-restriction:
given an arrangement ${\mathcal A}$ and a choice of hyperplane
$H \in {\mathcal A}$, set
$${\mathcal A'} = {\mathcal A} \setminus H \mbox{ and }{\mathcal A''} =
{\mathcal A}|_H.$$ The collection
$({\mathcal A}', {\mathcal A}, {\mathcal A}'')$ is called a {\em triple},
and a triple yields (see Proposition 4.45 of \cite{ot}) a left exact sequence
$$ 0 \longrightarrow D({\mathcal A}')(-1) \stackrel {\cdot H}{\longrightarrow}
D({\mathcal A}) \longrightarrow D({\mathcal A}'').$$ For a triple with ${\mathcal A} \subseteq \mathbb{P}^2$,
more is true (see \cite{cmh}): after pruning the Euler derivations and sheafifying, there is an exact sequence
\begin{equation}\label{eq:ses}
0 \longrightarrow {\mathcal D_0}'(-1) \longrightarrow {\mathcal D_0}
  \longrightarrow i_*{\mathcal D_0}'' \longrightarrow 0,
\end{equation}
where $i: H \hookrightarrow \mathbb{P}^2$; $i_*{\mathcal D_0}''\simeq {\mathcal O}_H(1-|{\mathcal A}''|).$
In \cite{t2}, Terao showed that freeness of a triple is related via:
\begin{thm}\label{thm:teraoAD}$($Addition-Deletion$)$
Let  $({\mathcal A}', {\mathcal A}, {\mathcal A}'')$ be a triple. Then any two of the following imply the third
\begin{itemize}
\item $D({\mathcal A})\simeq \oplus_{i=1}^n S(-b_i)$
\item $D({\mathcal A}')\simeq S(-b_n +1)\oplus_{i=1}^{n-1} S(-b_i)$
\item $D({\mathcal A}'')\simeq \oplus_{i=1}^{n-1} S(-b_i)$
\end{itemize}
\end{thm}
Theorem \ref{thm:teraoAD} applies in general, not just to arrangements in $\p^2$.
A smooth conic is intrinsically a $\mathbb{P}^1$, so it is natural
to ask if CL arrangements which admit a short exact sequence similar
to (\ref{eq:ses}) have an addition-deletion theorem; we tackle this
in the next two sections.

\pagebreak
A second criterion for freeness is special to the
case of line arrangements; to state it we need to define freeness for
{\em multiarrangements}. A multiarrangement $(\mathcal{A},{\bf m})$
is an arrangement together with a multiplicity $m_i$ for each
hyperplane. The module of derivations consists of $\theta$ such
that $l_i^{m_i} | \theta(l_i)$. As shown by Ziegler
in \cite{z}, freeness of multiarrangements is not combinatorial;
for recent progress see \cite{wy}.
\begin{thm}\label{thm:yoshi}$($Yoshinaga's multiarrangement criterion \cite{y1}$)$
$\mathcal{A} \subseteq \mathbb{P}^2$ is free iff
$\pi({\mathcal A},t) = (1+t)(1+at)(1+bt)$  and
$\forall H \in {\mathcal A}$ the multiarrangement
${\mathcal A}|_H$ has minimal generators in degree $a$ and $b$.
\end{thm}

The main results of this paper (Theorems \ref{thm:ADline}
and \ref{thm:ADconic}) show that
an addition-deletion construction holds for CL arrangements
with quasihomogeneous singularities; the
freeness of Example \ref{firstex} is explained by our results.
As one application, we show that a free CL arrangement, when
restricted to different lines, can yield multiarrangements
with different exponents; hence any version of Theorem \ref{thm:yoshi} for
CL arrangements will be quite subtle.
An addition-deletion theorem for multiarrangements
has recently been proven by Abe-Terao-Wakefield in \cite{atw};
our results are the first (to our knowledge) to give an inductive
criterion for freeness for nonlinear arrangements.


\section{Addition-Deletion for a line}
Let $(\mathcal C',\mathcal C, \mathcal C'')$ be a triple of CL
arrangements in $\mathbb P^2$, where $\mathcal C'=\mathcal C \setminus
\{L\}$,
$\mathcal C''=\mathcal C|_L$ and $L\in \mathcal C$ {\em is a line}. We begin by examining some examples:

\begin{exm} Let $\mathcal C'$ be the union of: $$\begin{array}{*{3}c}
C_1& = &x^2-xz+5y^2-5yz=0 \\
C_2& = &x^2+2y^2-xz-2yz=0 \\
L_1& = &x=0 \\
L_2& = &y=0 \\
L_3& = &x+y-z=0
\end{array}$$

\begin{figure}[h]
\begin{center}
\includegraphics{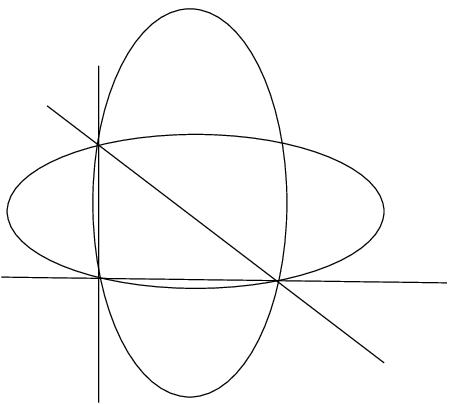}
\end{center}
\end{figure}

$D(\mathcal C')$ is free with exponents $\{1,2,4\}$, and the
degree of the Jacobian ideal is $28$, which is  equal to the
sum of the Milnor numbers at the intersection points.
Therefore at each singular point $\tau=\mu$.
If we restrict to any line, the corresponding multiarrangement
has two points of multiplicity 3, and it follows from \cite{wy} that
the exponents are $\{3,3\}$. Hence the obvious generalization of
Yoshinaga's criterion does not hold.
\end{exm}

\begin{exm}  Let $L_4=\{x-y=0\}$ and let ${\mathcal C_1} = {\mathcal C'} \cup L_4$.
\begin{figure}[h]
\begin{center}
\epsfig{file=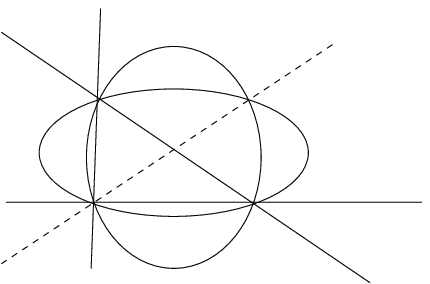,height=2.0in,width=2.0in}
\end{center}
\end{figure}
The degree of the Jacobian ideal is $39$, which is equal to the sum of Milnor
numbers at the points. It will follow from our results
that $D(\mathcal C_1)$ is free with
exponents $\{1,2,5\}$.
\end{exm}

\begin{exm} Let $L_4=\{x-2y=0\}$ and let ${\mathcal C_2} = {\mathcal
    C'} \cup L_4$. Then $\mathcal C_2$
  is free with exponents $\{1,3,4\}$. The degree of the Jacobian ideal
  is $37$, whereas the sum of the Milnor numbers is 38; the singularity at
$(0:0:1)$ has $\tau=15$ and $\mu=16$.
\begin{figure}[h]
\begin{center}
\epsfig{file=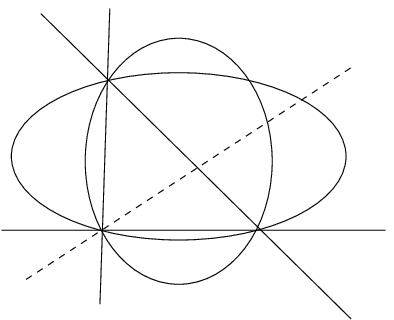}
\end{center}
\end{figure}
\end{exm}
\noindent For CL arrangements similar to ${\mathcal C_1}$, there is an addition-deletion theorem:
\begin{defn} A triple $(\mathcal C',\mathcal C, \mathcal C'')$ of CL
  arrangements is called quasihomogeneous if $\tau=\mu$ at each singular point of $\mathcal C'$ and $\mathcal C$.
\end{defn}

\begin{thm}\label{thm:ADline} Let $(\mathcal C',\mathcal C, \mathcal
  C'')$ be a quasihomogeneous triple with
$|L \cap \mathcal C|= |\mathcal C''| = k$.
The following are equivalent:
\begin{enumerate}
    \item $\mathcal C'$ is free with exponents $\{1,k-1,a\}$.
    \item $\mathcal C$ is free with exponents $\{1,k-1,a+1\}$.
\end{enumerate}
\end{thm}
\noindent Examples 2.1 and 2.2 illustrate the theorem; before giving the
proof of Theorem \ref{thm:ADline},
we need some preliminaries.
\begin{lem}\label{lem:derexact}
Let $L=\{x=0\}$. Then the maps
$p:D(\mathcal C')\longrightarrow D(\mathcal C)$, $p(\theta)=x\theta$
and $q:D(\mathcal C)\longrightarrow D(\mathcal C'')$,
$q(a\partial_x+b\partial_y+c\partial_z)=b(0,y,z)\partial_y+c(0,y,z)\partial_z$
are well defined and yield an exact sequence: $$0\longrightarrow D(\mathcal C')\longrightarrow D(\mathcal C)\longrightarrow D(\mathcal C'').$$
\end{lem}
\begin{proof} Let $f=xf'$ be the defining polynomial of $\mathcal C$, where $f'$ is the defining polynomial of
$\mathcal C'$. Then the defining polynomial of $\mathcal C''$ is
  $f''=Rad(f'|_{x=0})$. If $\theta'\in D(\mathcal C')$, then
  $\theta'(f')=Pf'$ for some $P\in S$;
\[
p(\theta')(f)=x\theta'(xf')= x(f'\theta'(x)+x\theta'(f'))\in\langle
f\rangle.
\]
So $p$ is well defined and injective.
Let $\theta=a\partial_x+b\partial_y+c\partial_z \in D(\mathcal C)$.
Then $\theta(x)=a\in \langle x \rangle$, so $a=xa'$.
If $\theta \in \ker(q)$, then $b=xb'$ and $c=xc'$, hence
$\theta=x\theta'$,
where $\theta'=a'\partial_x+b'\partial_y+c'\partial_z$.
Because $\theta\in D(\mathcal C)$,  $\theta(f')=x\theta'(f')
\in \langle f' \rangle$. Since $x$ and $f'$ are relatively prime, we
get that $\theta'(f') \in \langle f' \rangle$,
which implies that $\theta\in \im(p)$.

It remains to show is that $q$ is well defined. For suitable
$u_i,v_i\in \mathbb C$ and $m_i \in {\mathbb Z}$ we have that
\[
f'|_{x=0}=\prod_i(u_iy+v_iz)^{m_i}, \mbox{ so }
f''=\prod_i(u_iy+v_iz).
\]
Let $L'$ be a line in $\mathcal C'$ defined by the vanishing of
$t_ix+u_iy+v_iz=0$ for some $i$ and $t_i\in \mathbb C$, and let
$\theta=a\partial_x+b\partial_y+c\partial_z \in D(\mathcal C)$.
Then $\theta(L')\in \langle L' \rangle$, so evaluating at $x=0$
and using the earlier observation that $a=xa'$, we find
$(b(0,y,z)\partial_y+c(0,y,z)\partial_z)(u_iy+v_iz)\in \langle
u_iy+v_iz\rangle$.

Now suppose $C$ is a conic in $\mathcal C'$; after a change of
coordinates we may assume $C$ intersects $L=\{x=0\}$ in the points
$(0:0:1)$ and $(0:u:v)$. Then $C=xA+y(vy-uz)$ and
$C|_{x=0}=y(vy-uz)$, where $A$ is some linear form.
We have
\[
\theta(C)\!=\!a(A\!+\!x\partial_x(A))\!+\!x(b\partial_y(A)\!+\!c\partial_z(A))\!+\!b\partial_y(y(vy\!-\!uz))\!+\!c\partial_z(y(vy\!-\!uz))\!\in\!
\langle C \rangle.
\]
Evaluating at $x=0$ and again using that $a=xa'$ we find
\[
(b(0,y,z)\partial_y+c(0,y,z)\partial_z)(y(vy-uz))\in \langle
y(vy-uz)\rangle.
\]
Since $y$ and $vy-uz$ are relatively prime we obtain
\begin{center}
$\begin{array}{c}
(b(0,y,z)\partial_y+c(0,y,z)\partial_z)(y)\in \langle y\rangle \\
(b(0,y,z)\partial_y+c(0,y,z)\partial_z)(vy-uz)\in \langle
vy-uz\rangle.
\end{array}$
\end{center}
This shows that for each factor $u_iy+v_iz$ of $f''$,
\[
(b(0,y,z)\partial_y+c(0,y,z)\partial_z)(u_iy+v_iz)\in \langle
u_iy+v_iz\rangle,
\]
so the map $q$ is well defined. It follows that
$D_0(C'')=\mathbb C[y,z](-(k-1))$, where $k = |L \cap {\mathcal C}'| =
\mbox{deg}(f'')$. A similar argument works if $C$ is tangent to $L$.
\end{proof}

\begin{lem}\label{lem:ADmilnors}
Let $X$ and $Y$ be two reduced plane curves with no common component,
meeting at a point $p$.
Then
\[
\mu_p(X \cup Y) = \mu_p(X) + \mu_p(Y) +2(X \cdot Y)_p -1,
\]
where $(X \cdot Y)_p$ is the intersection number of $X$ and $Y$ at $p$.
\end{lem}
\begin{proof}
See \cite{w}, Theorem 6.5.1; the point is that the Milnor fiber is
a connected curve, and the result follows from using the additivity of
the Euler characteristic and the interpretation of $\mu_p$ as the
first betti number of the Milnor fiber.
\end{proof}
\begin{prop}\label{prop:sesLine}
 Let $(\mathcal C',\mathcal C, \mathcal C'')$ be a
 quasihomogeneous triple. Then
\[
 0\longrightarrow {\mathcal D_0}'(-1) \longrightarrow {\mathcal D_0}
 \longrightarrow i_*{\mathcal D_0}''.
\]
is also right exact.
\end{prop}

\begin{proof} It follows from Lemma \ref{lem:derexact} that
quotienting by the Euler derivation and sheafifying yields the
left exact sequence above; so it will suffice to show that
$HP(D_0,t)=HP(D_0'(-1),t)+HP(i_*D_0'',t)$, where $HP(-,t)$
denotes the Hilbert polynomial. For an CL
arrangement $\mathcal C$ with $m$
lines and $n$ conics, let $d =2n+m-1$. We have an exact sequence:
$$0\longrightarrow D_0(\mathcal C)\longrightarrow S^3\longrightarrow S(d)\longrightarrow S(d)/J
\longrightarrow 0,$$ where $S=\mathbb K[x,y,z]$ and $J$ is the
Jacobian ideal of the defining polynomial of $\mathcal C$. Since
\[
HP(D_0,t) = 3{{t+2}\choose{2}}-{{t+2+d}\choose{2}}+\deg(J)
\]
\[
HP(D_0'(-1),t) = 3{{t+1}\choose{2}}-{{t+d}\choose{2}}+\deg(J'),
\]
we find that
\[
HP(D_0,t)-HP(D_0'(-1),t)=\deg(J)-\deg(J')+t-2d+2.
\]
By the assumption that $(\mathcal C',\mathcal C, \mathcal C'')$ is a
quasihomogeneous triple,
\[
\deg(J) =\!\!\!\! \sum\limits_{p \in \Sing(C)}\!\!\!\! \mu_p(C) \mbox{ and }
\deg(J') =\!\!\!\! \sum\limits_{p \in \Sing(C')}\!\!\!\! \mu_{p}(C').
\]
Let $\alpha$ be the sum of Milnor numbers of points off $L$,
so
\[
\deg(J)=\alpha+\sum_{p \in L \cap C'}\mu_p(C).
\]
Since $\mu_p(L) =0$, by Lemma \ref{lem:ADmilnors}, the above is
\[
\alpha+\sum_{p \in L \cap C'}(\mu_p(C')+ 2(L \cdot C')_p -1).
\]
As $\deg(J')=\alpha+\sum_{p \in L \cap C'} \mu_p(C')$ and $\mid L \cap
C'\mid = k$, we obtain:
\[
\deg(J)-\deg(J')=2\!\!\sum_{p \in L \cap C'}\!\!(L \cdot C')_p -k.
\]
By  Bezout's theorem,
\[
\sum_{p \in L \cap C'}(L \cdot C')=d,
\]
so $\deg(J)-\deg(J')=2d-k$, hence
\[
HP(D_0,t)-HP(D_0'(-1),t)=t+2-k=t+1-(|\mathcal C''|-1).
\]
Since $i_*{\mathcal D_0}''=\mathcal O_L(1-|\mathcal C''|)$, this
yields the result.
\end{proof}
\begin{defn} A coherent sheaf $\mathcal F$ on $\mathbb P^r$ is $m-$regular iff
$H^i\mathcal F(m-i)=0$ for every $i\geq 1$. The smallest number $m$
  such that $\mathcal F$ is $m-$regular is $\reg(\mathcal F)$.
\end{defn}
\begin{lem}\label{lem:regbounds}
For a quasihomogeneous triple with $|\mathcal C''|=k$,
\[
\reg(\mathcal D_0)\leq \max\{\reg(\mathcal D_0')+1, k-1\}.
\]
\end{lem}
\begin{proof} Immediate from Proposition \ref{prop:sesLine} (see \cite{cmh}).
\end{proof}

\begin{lem}\label{lem:Lineexact}
If $\mathcal D_0'= \mathcal O_{\mathbb P^2}(1-k)\oplus \mathcal
O_{\mathbb P^2}(-a)$, then there is an exact sequence of $S$-modules:
\[
0\longrightarrow D_0'(-1)\longrightarrow D_0\longrightarrow
D_0''\longrightarrow 0.
\]
\end{lem}
\begin{proof}
For all $t$, $H^1(\mathcal D_0'(t-1))=0$, so the long exact sequence
in cohomology arising from Proposition \ref{prop:sesLine} and the vanishing of
$H^1(\mathcal D_0'(t))$ yield an exact sequence:
\[
0\longrightarrow \bigoplus_tH^0(\mathcal D_0'(-1)(t)) \longrightarrow \bigoplus_tH^0(\mathcal D_0(t))
\longrightarrow \bigoplus_tH^0(\mathcal D_0''(t))\longrightarrow 0.
\]
Theorem A.4.1 of \cite{e} relates a graded module to its sheaf and
local cohomology (at the maximal ideal $\mathfrak{m}$) modules:
\[
0\longrightarrow H_{\mathfrak{m}}^0(D_0)\longrightarrow D_0
\longrightarrow \bigoplus_tH^0(\mathcal D_0(t))\longrightarrow
H_{\mathfrak{m}}^1(D_0)\longrightarrow 0.
\]
This is true also for $D_0'(-1)$ and $D_0''$. By \cite{e}, A.4.3,
$H_\mathfrak{m}^0(M)=H_\mathfrak{m}^1(M)=0$ if $depth(M)\geq 2$. Lemma
2.1 of \cite{ms} gives the desired bound on depth for the modules of
derivations, which concludes the proof.
\end{proof}

\noindent The next two lemmas prove the two implications
$(1)\rightarrow (2)$ and $(2)\rightarrow (1)$ of Theorem \ref{thm:ADline}.
In what follows, $(\mathcal C',\mathcal C, \mathcal C'')$ is a
quasihomogeneous triple, with $L$ a line and $|L \cap \mathcal C|= |\mathcal C''| =
k$.

\begin{lem}\label{lem:Lines1}
If $\mathcal C'$ is free with $exp(\mathcal C') = \{1,k-1,a\}$, then
$\mathcal C$ is free with $exp(\mathcal C) = \{1,k-1,a+1\}$.
\end{lem}
\begin{proof}
First, if $S = \K[x_1,\ldots,x_{\ell}]$, and $S(-i)$ is
a free graded $S$-module with generator in degree $i$, then the Hilbert
series satisfies
\[
HS(S(-i),t) = \sum\limits_{j \in \Z} \dim_{\K}\!S(-i)_j \cdot  t^j =
\frac{t^i}{(1-t)^{\ell}}.
\]
If $\mathcal C'$ is free with exponents $\{1,k-1,a\}$, then
$D'_0(-1) \simeq S(-k) \oplus S(-1-a)$. It follows from the proof
of Proposition \ref{prop:sesLine} that
$HS(D_0'',t) = \frac{t^{k-1}}{(1-t)^2}$.
So by Lemma \ref{lem:Lineexact} and the additivity of
Hilbert series on an exact sequence,
\[
HS(D_0) =  \frac{t^{k}+t^{a+1}}{(1-t)^3} + \frac{t^{k-1}}{(1-t)^2} = \frac{t^{a+1}+t^{k-1}}{(1-t)^3}.
\]
Since $D_0' \simeq S(-k+1) \oplus S(-a)$,
$\reg(D_0')=\max\{k-1,a\}$. By Lemma \ref{lem:regbounds}, if $a \geq k-1$, then
$\reg(D_0)\leq a+1$; and if $a\leq k-2$, then $\reg(D_0)\leq k$.
If $a\leq k-2$, then a free resolution for $D_0$ is of the
form :
\[
0 \longleftarrow D_0 \longleftarrow S(-k+1) \oplus S(-a-1) \oplus S(-b)^d  \longleftarrow S(-b)^d \longleftarrow 0.
\]
From regularity constraints, $b$ must be at most $k$. As this is a minimal
free resolution, and it is impossible to have a syzygy on a single generator,
the only situation which can actually arise occurs when $b=k$:
\[
0 \longleftarrow D_0 \longleftarrow S(-k+1) \oplus S(-a-1) \oplus S(-k)^d  \longleftarrow S(-k)^d
\longleftarrow 0.
\]
Let $t_1, t_2$ be two independent derivations in $D_0$ of degrees
$\deg(t_1)=a+1$ and $\deg(t_2)=k-1$; our computation of the Hilbert
series, combined with the fact that $\pdim(D_0)\leq 1$ means such
derivations must exist. Let $E,t_1',t_2'$ be a basis for $D'$
with $\deg(t_1')=a$ and $\deg(t_2')=k-1$, and $E$ the Euler derivation.

Now note that $t_1'\in D_0'\setminus D_0$, for otherwise in $D_0$
there would be an element of degree $a$. So $t_1'(x)\notin \langle x \rangle$.
Since $D\subset D'$, then $t_1=f_1E+xt_1'$ and $t_2=f_2E+ut_2'+ft_1'$,
where $u$ is a constant, $\deg(f)=k-1-a$, $\deg(f_1)=a$ and
$\deg(f_2)=k-2$. For a resolution as above, $gt_1=Lt_2$, where $L$
is a linear form and $\deg(g)=k-(a+1)$. Hence
\[
(gf_1-Lf_2)E+(gx-Lf)t_1'+(-Lu)t_2'=0,
\]
and since $E,t_1',t_2'$ is a basis we find that $u$ vanishes and
$gx=Lf$. But $(t_2-f_2E)(x)\in \langle x \rangle$ and
$t_1'(x)\notin \langle x \rangle$. Since $u=0$, $x$ must divide $f$,
and so $g=Lg'$ for some $g'$. Since $gt_1=Lt_2$, we obtain
$t_2=g't_1$, a contradiction.
If $a\geq k-1$, simply switch the roles of $a$ and $k$ above.
\end{proof}

\begin{lem}\label{lem:Lines2}
If $\mathcal C$ is free with $exp(\mathcal C) = \{1,k-1,a+1\}$, then
$\mathcal C'$ is free with $exp(\mathcal C') = \{1,k-1,a\}$.
\end{lem}
\begin{proof}
In order to obtain an appropriate vanishing, we need to
dualize. Apply $Hom(-,\mathcal O_{\mathbb P^2})$
to the exact sequence
\[
0 \longrightarrow {\mathcal
D_0}'(-1) \longrightarrow {\mathcal D_0} \longrightarrow i_*{\mathcal
  D_0}'' \longrightarrow 0.
\]
The vanishing of $Hom_{{\mathcal O}_{\mathbb P^2}}({\mathcal O}_{\mathbb P^1}(t),{\mathcal O}_{\mathbb P^2})$ and $Ext_{{\mathcal O}_{\mathbb P^2}}^1(\mathcal D_0, \mathcal O_{\mathbb P^2})$ yield an exact sequence:
\[
0 \longrightarrow \mathcal D_0^{\vee} \longrightarrow \mathcal
D_0'^{\vee}(1) \longrightarrow Ext_S^1(\mathcal O_L(1-k), \mathcal
O_{\mathbb P^2})\longrightarrow 0.
\]
The free ${\mathcal O}_{\mathbb P^2}$ resolution for $\mathcal O_L(1-k)$ is:
\[
0\longrightarrow O_{\mathbb P^2}(-k)\longrightarrow O_{\mathbb
P^2}(1-k)\longrightarrow \mathcal O_L(1-k)\longrightarrow 0,
\]
so $Ext_S^1(\mathcal O_L(1-k), \mathcal{O}_{\mathbb{P}^2}) \simeq \mathcal O_L(k)$.
Since $\mathcal D_0^{\vee} =\mathcal O_{\mathbb P^2}(k-1)\oplus \mathcal
O_{\mathbb P^2}(a+1)$, combining this with the long exact sequence
in cohomology yields a regularity bound
\[
\reg(\mathcal D_0'^{\vee})\leq \max\{\reg(\mathcal D_0^{\vee})+1, 1-k\},
\]
and the exact sequence of $S-$modules:
\[
0\longrightarrow D_0^{\vee}(-1)\longrightarrow
D_0'^{\vee}\longrightarrow S(k-1)/L\longrightarrow 0,
\]
with $D_0^{\vee}=S(k-1)\oplus S(a+1)$. So:
\[
HS(D_0'^{\vee}) = \frac{t^{-a}+t^{1-k}}{(1-t)^3}.
\]
An argument as in the proof of Lemma \ref{lem:Lines1} shows that
$D_0'^{\vee}=S(a)\oplus S(k-1)$, hence $D(\mathcal C')$ is free
with exponents $\{1,k-1,a\}$.
\end{proof}
\begin{cor}
A free CL arrangement, when restricted to a line, can yield
different multiarrangements.
\end{cor}
\begin{proof}
In Example 2.2, add the line $L = \{x-\alpha y+(\alpha-1)z=0\}$, where $\alpha
\not \in \{0,1,-5,-2,\infty\}$. Then $L$ passes through $(1:1:1)$, and
the choices for $\alpha$ ensure that $L$ is not tangent to any conic,
and misses all singularities save
$(1:1:1)$. The new arrangement is quasihomogeneous, and $L$
meets $\mathcal C_1$ in six points. By
Theorem \ref{thm:ADline}, the new arrangement is free with exponents $\{1,3,5\}$.

Restrict this new arrangement to the line $L_3 = \{x+y-z=0\}$.
After a change of coordinates, we obtain a multiarrangement with
defining polynomial $$x^3y^3(x-y)(\alpha x-y).$$ This is exactly
Ziegler's example from \cite{z}: $\alpha =-1$ gives exponents $\{3,5\}$, and
for $\alpha \neq -1$, the exponents are $\{4,4\}$.
\end{proof}

\section{Addition-Deletion for a conic}
\noindent Let $(\mathcal C',\mathcal C, \mathcal C'')$ be a triple of CL
arrangements in $\mathbb P^2$, where $C$ is a conic in $\mathcal C$,
and ${\mathcal C}'={\mathcal C} \setminus \{ C \}$, ${\mathcal C}''={\mathcal
C}'|_C$.  We begin with some examples.

\begin{exm} Suppose $\mathcal C$ is as in Example 2.2, so
$\mathcal C$ has quasihomogeneous singularities, and is
free with exponents $\{1,2,5\}$. If we delete one of the conics, the
resulting arrangement ${\mathcal C}'$ is free and quasihomogeneous,
with exponents $\{1,2,3\}$.
\end{exm}

\noindent When $k$ is odd, the situation is more complicated:
\begin{exm} Let $\mathcal C'$ be the braid arrangement $A_3 = V(xyz(x-z)(y-z)(x+y-z))$, and $\mathcal C = \mathcal C' \cup C$,
where the conic $C = V(xy+7xz+13yz)$.
\begin{figure}[h]
\begin{center}
\epsfig{file=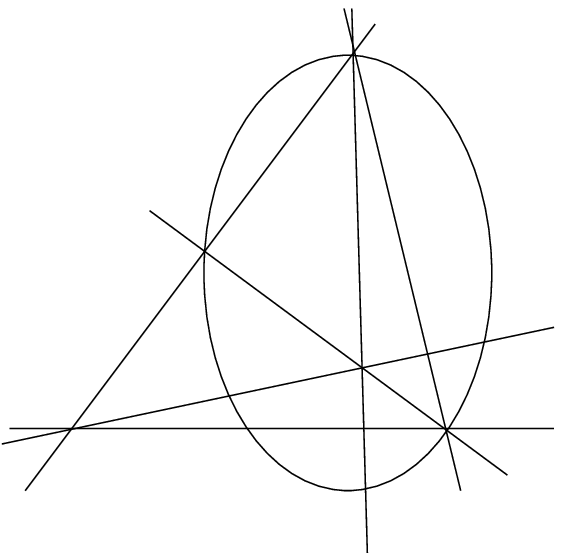,height=1.5in,width=1.8in}
\end{center}
\end{figure}

\noindent $\mathcal C'$ is a free arrangement with exponents
$\{1,2,3\}$, and $|\mathcal C''| = 7$. ${\mathcal C}$ is also
quasihomogeneous, but not free.
\end{exm}

\begin{exm} Let  $\mathcal C$ be the quasihomogeneous CL arrangement
with defining polynomial $(x^2-xz+2y^2-2yz)xy(x+y-z)$.
\begin{figure}[h]
\begin{center}
\epsfig{file=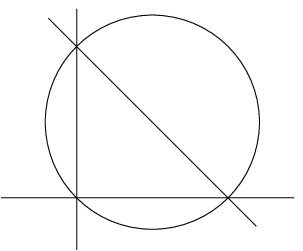,height=1.3in,width=1.6in}
\end{center}
\end{figure}
$D(\mathcal C)$ is free with exponents $\{1,2,2\}$. Deleting the conic
yields a free line arrangement with exponents $\{1,1,1\}$.
\end{exm}
\begin{thm}\label{thm:ADconic}
Let $(\mathcal C',\mathcal C, \mathcal C'')$ be a quasihomogeneous
triple, with $|C \cap \mathcal C'|\! =\!|\mathcal C''|\! =\! k$.
If $k=2m$ then the following are equivalent:
\begin{enumerate}
\item $\mathcal C'$ is free with $exp(\mathcal C') = \{1,m,a\}$.
\item $\mathcal C$ is free with $exp(\mathcal C) = \{1,m,a+2\}$.
\end{enumerate}
If $k=2m+1$ then:
\begin{enumerate}
\item $exp(\mathcal C') = \{1,m,m \}  \Longleftrightarrow exp(\mathcal C)
= \{1,m+1,m+1 \} $.
\item if $exp(\mathcal C') = \{1,m,a \}$ with $a \ne m$ then $\mathcal C$
is not free.
\item if $exp(\mathcal C) = \{1,m+1,a+1 \}$ with $a \ne m$ then $\mathcal C'$
is not free.
\end{enumerate}
\end{thm}

We begin with some preliminaries. After an appropriate change of
coordinates, we may suppose that $C=\{y^2-xz=0\}$. Let $i$ be the composition of the maps $$\mathbb P^1\stackrel{v}{\longrightarrow} C
\hookrightarrow \mathbb P^2,$$ where $v(s:t)=(s^2:st:t^2)$, and
let $\psi$ be the composite map:
$$S={\mathbb K}[x,y,z] \stackrel{\phi}{\longrightarrow}\mathbb K[s^2,st,t^2] \hookrightarrow \mathbb K[s,t],$$ where $\phi(x)=s^2,
\phi(y)=st, \phi(z)=t^2.$

\vskip .1in

Let $\theta = a_1\partial_x+a_2\partial_y+a_3\partial_z\in D(\mathcal
C)$ be a derivation. Then $\theta(C)\in \langle C \rangle$, which
means $-za_1+2ya_2-xa_3=(y^2-xz)P$ for some $P\in S$. Via the map $\psi$ this
translates into $$t^2\psi(a_1)-2st\psi(a_2)+s^2\psi(a_3)=0.$$ So there exist $Q_1,Q_2\in \mathbb K[s,t]$ such
that $$\begin{array}{*{3}c}
\psi(a_1)& = &sQ_1 \\
\psi(a_2)& = &\frac{tQ_1+sQ_2}{2} \\
\psi(a_3)& = &tQ_2.
\end{array}$$

\vskip .1in

If $\psi:S\longrightarrow A$ is a ring map and $M$ is an $A-$module,
let $M_{\psi}$ denote the $S-$module obtained by restriction of
scalars.

\begin{prop} There is an exact sequence of $S-$modules $$0\longrightarrow D(\mathcal C')(-2) \stackrel{\cdot
C}{\longrightarrow} D(\mathcal C)\stackrel{\rho}{\longrightarrow} D(\mathcal A'')_{\psi},$$ where
$$\rho(a_1\partial_x+a_2\partial_y+a_3\partial_z)=Q_1\partial_s+Q_2\partial_t,$$ for every $a_1\partial_x+a_2
\partial_y+a_3\partial_z\in D(\mathcal C)$ and $Q_1,Q_2$ are defined as above; and $\mathcal A''$ is the
arrangement of the reduced points $i^{-1} (C\cap \mathcal C')$ in $\mathbb P^1$.
\end{prop}

\begin{proof} It is easy to check that $\rho$ is a
  homomorphism. For exactness, note:
\begin{center}
$\begin{array}{ccc}
\theta = a_1\partial_x+a_2\partial_y+a_3\partial_z\in \ker(\rho) &
  \leftrightarrow & Q_1=0 \mbox{ and }Q_2=0 \\
  & \leftrightarrow & \psi(a_1)=\psi(a_2)=\psi(a_3)=0 \\
& \leftrightarrow & a_1,a_2,a_3\in\langle y^2-xz\rangle \\
& \leftrightarrow & \theta=C \theta' \mbox{ with } \theta' \in D(\mathcal C').
\end{array}$
\end{center}
It remains to show that the image of $\rho$ is in $D(\mathcal A'')_{\psi}$.
Suppose $\alpha x+\beta y+\gamma z=0$ is a line of $\mathcal C$.
Let $\theta =a_1\partial_x+a_2\partial_y+a_3\partial_z\in D(\mathcal C)$. Then $$\alpha a_1+\beta a_2+\gamma a_3 = (\alpha
x+\beta y+\gamma z)P_1$$ for some $P_1\in S$. Therefore
\[
\alpha \psi(a_1)+\beta \psi(a_2)+\gamma \psi(a_3) = (\alpha s^2+\beta
st+\gamma t^2)\psi(P_1),
\]
which implies
\[
(2\alpha s+\beta t)Q_1 +(\beta
s+2\gamma t)Q_2 = 2(\alpha s^2+\beta st+\gamma t^2)\psi(P_1).
\]
This means that $(Q_1\partial_s+Q_2\partial_t)(\alpha s^2+\beta st+\gamma t^2)\in (\alpha s^2+\beta
st+\gamma t^2)\mathbb K[s,t]$. Since $\alpha s^2+\beta st+\gamma t^2$ is the defining polynomial of the two points
$i^{-1}(\{\alpha x+\beta y+\gamma z=0\}\cap C)$ in $\mathbb P^1$, we get that $Q_1\partial_s+Q_2\partial_t$ is a
derivation on the arrangement of these two points.

\vskip .1in

Suppose $C'=\{u_0x^2+u_1xy+u_2xz+u_3y^2+u_4yz+u_5z^2=0\}$ is a conic in the CL arrangement $\mathcal C$. Let
$\theta = a_1\partial_x+a_2\partial_y+a_3\partial_z\in D(\mathcal
C)$. Computations as above show that
\[
(Q_1\partial_s+Q_2\partial_t)(u_0s^4+u_1s^3t+(u_2+u_3)s^2t^2+u_4st^3+u_5t^4)\]
\[
\in
(u_0s^4+u_1s^3t+(u_2+u_3)s^2t^2+u_4st^3+u_5t^4)\mathbb K[s,t].
\]
Since
$u_0s^4+u_1s^3t+(u_2+u_3)s^2t^2+u_4st^3+u_5t^4$ is the defining
polynomial of the four points $i^{-1}(C'\cap C)$ in $\mathbb P^1$, we
get that $Q_1\partial_s+Q_2\partial_t$ is a derivation on the
arrangement of these four points. Similar arguments work in the case
of tangencies.
\end{proof}

Let $\theta = a_1\partial_x+a_2\partial_y+a_3\partial_z\in D(\mathcal C)_d$ such that $\rho(\theta)=s\partial_s
+ t\partial_t$. Then $a_1,a_2,a_3\in S_d$ with
$\psi(a_1)=s^2,\psi(a_2)=st,\psi(a_3)=t^2$. Thus $d=1$ and
$\theta$ is the Euler derivation in $D(\mathcal C)$. So quotienting by
the Euler derivations yields an exact sequence:
\[
0\longrightarrow D_0'(-2)
\stackrel{\cdot C}{\longrightarrow}
D_0\stackrel{\rho}{\longrightarrow} (D(\mathcal A'')_0)_{\psi}.
\]
Since $|\mathcal A''|=k$, after sheafifying,
${\mathcal D}_0({\mathcal A}'')=\mathcal O_{\mathbb P^1}(-k)$, and
  hence the sheafification of $D_0(\mathcal A'')_{\psi}$ is
$i_*\mathcal O_{\mathbb P^1}(-k)$.

\begin{lem}\label{lem:HPconic} $HP(i_*\mathcal O_{\mathbb P^1}(-k),t)=2t+1-k$.
\end{lem}
\begin{proof}
\noindent CASE 1: $k=2m$. Let $E$ be the divisor of the reduced $k$ points $i^{-1}(C\cap \mathcal C')$. Then the ideal
sheaf $\mathcal I_E=\langle f\rangle$, where $f\in\mathbb K[s,t]$ of
degree $k=2m$. There exists $g \in S_m$, unique modulo
$(y^2-xz)$, such that $g(s^2,st,t^2)=f$.
Clearly $y^2-xz$ cannot divide $g$, otherwise $g(s^2,st,t^2)=0=f$, so
the ideal of the reduced $k$ points on $C$ is $\langle
y^2-xz,g\rangle$.
Hence $i_*\mathcal I_E=\langle \bar{g} \rangle$ as an ideal of $S/\langle
y^2-xz\rangle$. As an $S-$module, it has free resolution
\[
0\longrightarrow S(-2-m)\stackrel{\cdot
C}{\longrightarrow} S(-m) \longrightarrow \langle \bar{g} \rangle
\longrightarrow 0,
\]
which yields:
$$HP(i_*\mathcal O_{\mathbb P^1}(-2m),t)={{t+2-m}\choose{2}}-{{t-m}\choose{2}}=2t+1-2m.$$

\vskip .1in

\noindent CASE 2: $k=2m+1$. Let $E$ be the divisor of the reduced $k$ points $i^{-1}(C\cap \mathcal C')$. Then the ideal
sheaf $\mathcal I_E=\langle f\rangle$, where $f\in\mathbb \mathbb
K[s,t]$ of degree $k=2m+1$. Let $L_1, L_2 \in \mathbb K[s,t]_1$ be two
independent linear forms which do not divide $f$, and let $f_i=L_if$.
Since $\langle L_1f, L_2f\rangle=\langle L_1,L_2\rangle\cap \langle
f\rangle$, then $\langle f_1,f_2\rangle$ defines the same ideal sheaf on $\mathbb P^1$ as $\langle f \rangle$. So
$\mathcal I_E=\langle f_1,f_2\rangle$.

Both $f_1$ and $f_2$ are of even degree $2m+2$. So there exist
$g_1,g_2\in S=\mathbb K[x,y,z]$ of degree $m+1$ such that
$g_i(s^2,st,t^2)=f_i, i=1,2$.
Next we show that $J = \langle y^2-xz,g_1,g_2\rangle$ is the ideal
of the reduced points $C\cap\mathcal C'$ on $C$. To see this,
note that if $p\in C\cap\mathcal C'$, then $f_i(i^{-1}(p))=0, i=1,2$.
So $g_i(p)=0, i=1,2$, and hence $g_i\in J, i=1,2$. Clearly
$y^2-xz$ does not divide $g_i$, otherwise $f_i$ is identically zero.
Also, suppose $g_2=\lambda g_1+P(y^2-xz)$, where $\lambda$ is a
constant.
Then $f_2=\lambda f_1$, i.e. $L_2=\lambda L_1$; a contradiction.
So $J$ is the ideal of $2m+1$ points on the conic $y^2-xz=0$. By the
Hilbert-Burch theorem, such an ideal is minimally generated by the
$2\times 2$ minors of
\begin{center}
$\left[
\begin{array}{c}
x\, y \\
y\, z \\
\alpha\, \beta
\end{array}
\right]$
\end{center}
where both $\alpha$ and $\beta$ have degree $m$.
So indeed $\langle y^2-xz,g_1,g_2\rangle = J$, and
$i_*\mathcal I_E=\langle \bar{g_1},\bar{g_2}\rangle \subseteq
S/\langle y^2-xz\rangle$. As an $S-$module it has free
resolution $$0\longrightarrow S^2(-2-m)\stackrel{\tiny{ \left[
\begin{array}{c}
x\, y \\
y\, z
\end{array}
\right]} }{\longrightarrow} S^2(-1-m)\longrightarrow \langle \bar{g_1},\bar{g_2} \rangle \longrightarrow 0,$$ so
for the odd case we find that
$$HP(i_*\mathcal O_{\mathbb P^1}(-2m-1),t)=2{{t+1-m}\choose{2}}-2{{t-m}\choose{2}}$$
$$=2t-2m=2t+1-(2m+1).$$
\end{proof}

\begin{prop}\label{prop:sheafconicExact}
For a quasihomogeneous triple  $(\mathcal C',\mathcal C, \mathcal
C'' = \mathcal C'|_C )$, the sequence
\[
0\longrightarrow \mathcal D_0'(-2)\stackrel{\tiny{\cdot C}}\longrightarrow \mathcal D_0\longrightarrow
i_*\mathcal O_L(-k)\longrightarrow 0
\]
is exact, where $i:L \stackrel{[s^2:st:t^2]}{\longrightarrow} {\mathbb P}^2$.
\end{prop}
\begin{proof}
We have $HP(D_0,t)-HP(D_0'(-2),t)=2t-4d+9+(\deg(J_f)-\deg(J_{f'}))$,
where $d+1$ is the degree of the defining
polynomial $f$ of $\mathcal C$ and $f'$ is the defining polynomial of
$\mathcal C'$. Since $(\mathcal C',\mathcal C, \mathcal C'')$ is
a quasihomogeneous triple, Bezout's theorem and Lemma
\ref{lem:ADmilnors}
imply that $\deg(J_f)-\deg(J_{f'})=4d-4-k$ and hence
$$HP(D_0,t)-HP(D_0'(-2),t)=2t+1-k.$$
By Lemma \ref{lem:HPconic}, this is exactly the Hilbert polynomial of the sheaf
$i_*{\mathcal O}_{{\mathbb P}^1}(-k)$ associated to
$D_0({\mathcal C}'')_{\psi}$.
\end{proof}

\begin{lem}\label{lem:moduleconicExact} For a quasihomogeneous triple such that
$\mathcal C'$ is free with exponents $\{1,m,a\}$,
\[
0 \longrightarrow D_0'(-2) \longrightarrow
D_0 \longrightarrow  D_0(\mathcal A'')_{\psi} \longrightarrow 0
\]
is exact.
\end{lem}
\begin{proof}
As we've seen, $\bigoplus_t H^0((i_*\mathcal O_L(-k))(t))\!=\!\bigoplus_t
H^0(\mathcal O_{\mathbb P^1}(2t-k))\!=\! D_0(\mathcal A'')_{\psi}$.
With the assumption on $\mathcal C'$, $H^1(\mathcal D_0'(t-2))$
vanishes for all $t$, and exactness follows as in the
proof of Lemma \ref{lem:Lineexact}.
\end{proof}
\noindent Theorem \ref{thm:ADconic} will follow from the next two lemmas.
\begin{lem}\label{lem:Conic1}
Let $(\mathcal C',\mathcal C, \mathcal C'')$ be a quasihomogeneous
triple, with $|C \cap \mathcal C'|\! =\!|\mathcal C''|\! =\! k$.
If $\mathcal C'$ is free with exponents $\{1,m,a\}$, then
\begin{enumerate}
\item If $k=2m$ then $\mathcal C$ is free with $exp(\mathcal C) = \{1,m,a+2\}$.
\item If $k=2m+1$ and $a=m$, then $\mathcal C$ is free with
  $exp(\mathcal C) = \{1,m+1,m+1\}$.
\item If $k=2m+1$ and $a\ne m$, then $\mathcal C$ is not free.
\end{enumerate}
\end{lem}
\begin{proof}
It follows from the computations in the proof of Lemma \ref{lem:HPconic} that
\begin{itemize}
    \item If $k=2m+1$, then $HS( D_0(\mathcal A'')_{\psi},t)=\frac{2t^{m+1}}{(1-t)^2}$.
    \item If $k=2m$, then $HS( D_0(\mathcal A'')_{\psi},t)=\frac{t^m(1+t)}{(1-t)^2}$.
\end{itemize}
Combining these results yields the Hilbert series of $D_0$.
\pagebreak

\noindent CASE 1: $k=2m$.
By Lemma \ref{lem:moduleconicExact},
\[
HS(D_0,t)=\frac{t^m+t^{a+2}}{(1-t)^3}.
\]
Since $\pdim(D_0)\leq 1$, this means that there exist minimal
generators $\theta,\eta \in D_0$ with $\deg(\theta)=m$ and
$\deg(\eta)=a+2$. Suppose $\{E,\theta_1,\theta_2\}$ basis
for $D'$ with $E$ the Euler derivation and $\deg(\theta_1)=m,
\deg(\theta_2)=a$. We now use that $D\subset D'$.

\begin{itemize}
    \item $m<a$. Since $\theta\in D_0$, $\theta=fE+c\theta_1$ for some
    $c \in \C^{*}$. Then $\{E,\theta,\theta_2\}$ is a basis for $D'$, so by Saito's criterion $\{E,\theta,C\theta_2\}$ is a basis for $D$.
    \item $m=a$. Then $\theta=fE+c_1\theta_1+c_2\theta_2$, where $c_1,c_2$ constants, not both zero. If $c_2\neq
    0$, then $\{E,\theta_1,\theta\}$ is a basis for $D'$, so by
    Saito's criterion $\{E,C\theta_1,\theta\}$ is a basis for $D$.
    \item $m=a+1$. Then $\theta=fE+c_1\theta_1+L_2\theta_2$, where $c$ is a constant and $L_2$ is a linear form,
    not both zero. If $c_1=0$ then $L_2\theta_2(C)\in \langle C \rangle$. Since $C$ is irreducible, then
    $\theta_2(C)\in \langle C \rangle$, and so $\theta_2\in D_0$ is of
    degree $a<m,a+2$. This is inconsistent with the Hilbert series of
    $D_0$. So $c_1\neq 0$, and so $\{E,\theta,\theta_2\}$ is a basis
    for $D'$, and again by
    Saito's criterion $\{E,\theta,C\theta_2\}$ is a basis for $D$.
    \item $m=a+2$. Then $\theta=f_1E+c_1\theta_1+g_1\theta_2$, where $c_1$ is a constant and $g_1$ is a quadratic form,
    not both zero and $\eta=f_2E+c_2\theta_1+g_2\theta_2$, where $c_2$ is a constant and $g_2$ is a quadratic form,
    not both zero. If $c_1=c_2=0$, then either $g_i=c_i'C, c_i'\neq
    0,i=1,2$ or $\theta_2\in D_0$, a contradiction (because
    $\deg(\theta_2)=a$). Therefore $c_2'\theta-c_1'\eta\in D_0\cap
    ES=\{0\}$, contradicting the fact that $\theta,\eta$ are minimal generators of $D_0$.
    So if $c_2\neq 0$, then $\{E,\eta,\theta_2\}$ is a basis for $D'$,
    and so by
    Saito's criterion $\{E,\eta,C\theta_2\}$ is a basis for $D$.
    \item $m>a+2$. Then $\theta=f_1E+c_1\theta_1+g_1\theta_2$, where $c_1$ is a constant and $g_1$ is a polynomial,
    not both zero and $\eta=f_2E+c_2\theta_1+g_2\theta_2$, where $c_2$ is a constant and $g_2$ is a quadratic form,
    not both zero. If $c_1=c_2=0$, then $g_1=Cg_1', g_1'\neq 0$ and
    $g_2=c_2'C$, $c_2'$ nonzero constant, and the argument used
    above yields a contradiction. So $c_1\neq 0$ or $c_2\neq
    0$. Applying Saito's criterion yields the desired result.
\end{itemize}

\vskip .1in

\noindent CASE 2: $k=2m+1, m=a$.
By Lemma \ref{lem:moduleconicExact},
\[
HS(D_0,t)=\frac{2t^{m+1}}{(1-t)^3}.
\]
This implies there exist degree $m+1$ minimal generators
$\theta,\mu \in D_0$. Suppose $\{E,\theta_1,\theta_2\}$
is a basis for $D'$ where $E$ is the Euler derivation and $\deg(\theta_1)=m,
\deg(\theta_2)=m$. So $\theta=f_1E+L_1\theta_1+K_1\theta_2$ and
$\mu=f_2E+L_2\theta_1+K_2\theta_2$, where $L_1,L_2,K_1,K_2$ are linear
forms, and for any $i=1,2$, $L_i,K_i$ cannot be simultaneously
zero. Hence
$L_2\theta-L_1\mu-(L_2f_1-L_1f_2)E=(L_2K_1-L_1K_2)\theta_2 \in
D(\mathcal C)$.
But $\theta_2$ is in $D(\mathcal C')$ and
$\theta_2(C)\notin\langle C\rangle$, else $(D_0)_m$ is nonzero,
which is inconsistent with the Hilbert Series.
Hence $L_2K_1-L_1K_2=cC$, where $c$ is a constant.

If $c=0$, then $L_1=uK_1, L_2=uK_2, u\neq 0$ or $L_1=vL_2,
K_1=vK_2,v\neq0$, where $u,v$ are constants, and that
$K_2f_1=K_1f_2$ and $L_2f_1=L_1f_2$. If $L_1=uK_1, L_2=uK_2, u\neq 0$, and $K_1\neq ct\cdot K_2$ we get
$\theta=K_1(gE+u\theta_1+\theta_2)$. Since $K_1\neq 0$ (else $L_1=0$) then $\theta(C)\in\langle
C\rangle$ implies $(gE+u\theta_1+\theta_2)(C)\in\langle
C\rangle$, yielding a degree $m$ derivation in $D(\mathcal
C)$, a contradiction. If $K_1=ct\cdot K_2$, then
$\theta$ and $\mu$ are not minimal generators, also a contradiction.
If $c\neq 0$, then we find $\det[E,\theta,\mu]=cC\det[E,\theta_1,\theta_2]$,
and Saito's criterion shows that $\{E,\theta,\mu\}$ is a
basis for $D(\mathcal C)$.
\vskip .1in
\noindent CASE 3: $k=2m+1, m \ne a$.
By Lemma \ref{lem:moduleconicExact},
\[
HS(D_0,t)=\frac{t^{a+2}+2t^{m+1}-t^{m+2}}{(1-t)^3}.
\]
Since $m \ne a$, there is no cancellation in the numerator, hence
$D_0$ cannot be free.
\end{proof}

\begin{lem}\label{lem:Conic2}
Let $(\mathcal C',\mathcal C, \mathcal C'')$ be a quasihomogeneous
triple, with $|C \cap \mathcal C'|\! =\!|\mathcal C''|\! =\! k$.
If $\mathcal C$ is free, then
\begin{enumerate}
\!\!\item If $k\!=\!2m$ and $exp(\mathcal C)\!=\!\{1,m,a\!+\!2\}$, then $\mathcal C'$ is free with $exp(\mathcal C')\!=\!\{1,m,a\}$.
\!\!\item If $k\!\!=\!\!2m\!+\!1$ and $exp(\mathcal
  C)\!\!=\!\!\{1,m\!\!+\!\!1,m\!\!+\!\!1\}$, then $\mathcal C'$ is free with $exp(\mathcal C')\!=\!\{1,m,m\}$.
\!\!\item If $k\!=\!2m\!+\!1$ and $exp(\mathcal C)\!=\!\{1,m\!+\!1,a\!+\!1\}$ with $a\!\ne\!m$, then $\mathcal C'$ is not free.
\end{enumerate}
\end{lem}
\begin{proof}
As in Lemma \ref{lem:Lines2}, apply
$Hom(-,\mathcal O_{\mathbb P^2})$ to the exact sequence
\[
0\longrightarrow \mathcal D_0'(-2)\stackrel{\tiny{\cdot C}}\longrightarrow \mathcal D_0\longrightarrow
i_*\mathcal O_L(-k)\longrightarrow 0.
\]
Since $i_*\mathcal O_L(-k)$ is supported on
the conic $C$, $Hom(i_*\mathcal O_L(-k),\mathcal O_{\mathbb P^2})=0$.
The assumption that $D_0$ is free implies that
$Ext_S^1(\mathcal D_0, \mathcal O_{\mathbb P^2})=0$. This yields an
exact sequence:
\[
0
\longrightarrow \mathcal D_0^{\vee} \longrightarrow \mathcal D_0'^{\vee}(2) \longrightarrow Ext_S^1(i_*\mathcal
O_L(-k), \mathcal O_{\mathbb P^2})\longrightarrow 0.
\]
As $D_0$ free with known exponents, so also is $D_0^{\vee}$, and the
Hilbert Series is known. The proof of Lemma \ref{lem:HPconic} provides a
free resolution of $i_*\mathcal I_E$, which allows us to compute
$Ext_S^1(i_*\mathcal I_E,S)$. Combining everything yields the
Hilbert Series of $D_0'^{\vee}$, and the result follows as in the
previous analysis.
\end{proof}
\section{Freeness of CL arrangements is not combinatorial}
We close with a pair of examples which show that in the CL case,
Terao's conjecture that freeness is a combinatorial invariant of
an arrangement is false.
\begin{exm}\label{exm:last}
Let $\mathcal C_1$ be given by
\[
\begin{array}{*{3}c}
C_1& = &\{y^2+xz=0\} \\
C_2& = &\{y^2+x^2+2xz=0\} \\
L_1& = &\{x=0\}
\end{array}
\]
$L_1$ is tangent to both $C_1$ and $C_2$ at the point $P=(0:0:1)$;
$C_1$ and $C_2$ have two other points in common.
Adding the line $L=\{y=0\}$ passing through $P$ to $\mathcal C_1$
yields a
quasihomogeneous, free CL arrangement $\mathcal C$,
with $D_0(\mathcal C)\simeq S(-2)\oplus S(-3):$

\begin{center}
\epsfig{file=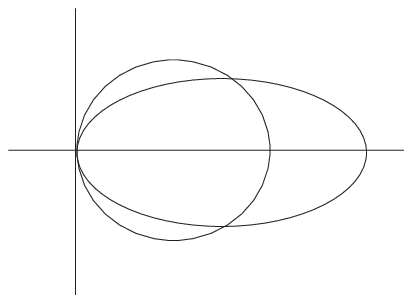,height=1.5in,width=1.8in}
\end{center}

\noindent The line $L'=\{x-13y=0\}$ passes through $P$, and
misses the other singularities of $\mathcal C_1$. The CL
arrangement $\mathcal C' = L' \cup \mathcal C_1$ is
combinatorially identical to $\mathcal C$, but $\mathcal C'$ is
not quasihomogeneous, and not free:
\[
0\longrightarrow S(-9)\longrightarrow S(-8)^3\longrightarrow
D_0(\mathcal C') \longrightarrow 0.
\]
\end{exm}
\noindent Even for CL arrangements with ordinary singularities,
there are counterexamples:
\begin{exm}\label{exm:noterao}
Let $\mathcal A$ be the union of the five smooth conics:
$$\begin{array}{*{3}c}
C_1& = &(x-3z)^2+(y-4z)^2-25z^2=0 \\
C_2& = &(x-4z)^2+(y-3z)^2-25z^2=0 \\
C_3& = &(x+3z)^2+(y-4z)^2-25z^2=0 \\
C_4& = &(x+4z)^2+(y-3z)^2-25z^2=0 \\
C_5& = &(x-5z)^2+y^2-25z^2=0
\end{array}$$
$\mathcal A$ has 13 singular points, all ordinary. At 10 of these
points only two branches of $\mathcal A$ meet, while at the
points $(0:0:1),(1:i:0),(1:-i:0)$, all five conics meet.
The Milnor and Tjurina numbers agree at all singularities
except $(0:0:1)$, where $\tau=15$ and $\mu=16$.
Adding lines $L_1,L_2,L_3$ connecting $(0:0:1),(1:i:0),(1:-i:0)$
yields a free CL arrangement $\mathcal C$, with $D_0(\mathcal
C) = S(-6)^2$.
\begin{figure}[h]
\begin{center}
\epsfig{file=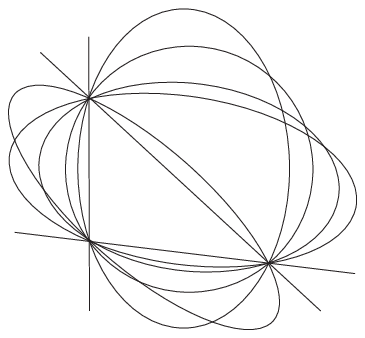,height=2.0in,width=2.0in}
\end{center}
\end{figure}
\vskip .05in
\noindent Next, let $\mathcal A'$ be the union of the following five smooth conics:
$$\begin{array}{*{3}c}
C_1& = &x^2+8y^2+21xy-xz-8yz=0 \\
C_2& = &x^2+5y^2+13xy-xz-5yz=0 \\
C_3& = &x^2+9y^2-4xy-xz-9yz=0 \\
C_4& = &x^2+11y^2+xy-xz-11yz=0 \\
C_5& = &x^2+17y^2-5xy-xz-17yz=0
\end{array}$$
$\mathcal A'$ is combinatorially identical to $\mathcal A$, but
at the points $(0:0:1),(1:0:1),(0:1:1)$ where all
the branches meet, $\tau=15$ and $\mu=16$.
Adding the lines
connecting these three points yields an CL arrangement $\mathcal
C'$ which is combinatorially identical to $\mathcal C$ but
{\em not free}; the free resolution of $D_0(\mathcal C')$ is:
\[
0\longrightarrow S(-8)^2 \longrightarrow S(-7)^4 \longrightarrow D_0(\mathcal
C') \longrightarrow 0.
\]
As was pointed out by the referee, the complements of
arrangements $\mathcal{A}$ and  $\mathcal{A}'$ are homeomorphic
(via a Cremona transformation centered on the three multiple
intersection points) to the complements of a pair of line arrangements
consisting of eight lines in general position. The moduli space of
such objects is connected, so the complements are rigidly isotopic,
hence homeomorphic. So freeness is also not a topological invariant.
\end{exm}
\pagebreak
\noindent {\bf Concluding Remarks}:
\begin{enumerate}
\item As noted in \S1.2, for the complement $X$ of an CL arrangement
in ${\mathbb P}^2$ the betti numbers $h^1(X)$ and $h^2(X)$ depend
only on the combinatorics, and so if $X$ is quasihomogeneous
and free, there is a version of Terao's theorem, which we leave for
the interested reader.

\item In the examples above, the Jacobian
ideals are of different degrees,
so are not even members of the same Hilbert scheme.
Do there exist CL arrangements with isomorphic intersection
poset {\em and} singularities which are locally isomorphic,
one free and one nonfree? Do there exist counterexamples
where all singularities are quasihomogeneous?

\item As shown by Example 2.3, quasihomogenity is not a necessary
condition for freeness of CL arrangements. However, without
this assumption, the sequences in Propositions 2.8 and 3.8 may
not be exact, which means that any form of addition-deletion will
require hypotheses on higher cohomology.
\end{enumerate}
\noindent {\bf Acknowledgments}:  Macaulay2 computations were
essential to our work. We also thank an anonymous referee for
many useful suggestions, in particular for pointing out that
we should remove one of our original conditions (that
$\mathcal{C}$ has only ordinary singularities).
\renewcommand{\baselinestretch}{1.0}
\small\normalsize 

\bibliographystyle{amsalpha}

\end{document}